\documentclass[a4paper,leqno,11pt]{article}
\usepackage{amsmath, amssymb, amsthm, amsfonts}
\usepackage{mathrsfs}

\theoremstyle{plain}
\newtheorem{theorem}{Theorem}

\newtheorem{lemma}{Lemma}

\theoremstyle{remark}
\newtheorem{remark}{Remark}

\renewcommand{\Re}{{\rm Re\,}}

\newcommand{\tlg}{t \log\frac{t}{2\pi}-t}

\numberwithin{equation}{section}

\title{Some mean value results related to Hardy's function}
\author{Xiaodong Cao, Yoshio Tanigawa and Wenguang Zhai\thanks{This work is  
supported by the National Key Basic Research Program of China(Grant No. 2013CB834201)}}
\date{\empty}

\begin{document}

\maketitle
\footnote[0]{2000 Mathematics Subject Classification: 11M06, 11N37.} 
\footnote[0]{Key words and phrases: Hardy's function, mean value theorems, approximate functional equation, exponential 
sum and integral}

\begin{abstract}
Let $\zeta(s)$ and $Z(t)$ be the Riemann zeta function and Hardy's function respectively. We show asymptotic formulas 
for $\int_0^T Z(t)\zeta(1/2+it)dt$ and $\int_0^T Z^2(t) \zeta(1/2+it)dt$. 
Furthermore we derive an upper bound for $\int_0^T Z^3(t)\chi^{\alpha}(1/2+it)dt$ for 
$-1/2<\alpha<1/2$, where $\chi(s)$ is the function which appears in the functional equation of 
the Riemann zeta function: $\zeta(s)=\chi(s)\zeta(1-s)$. 
\end{abstract}

\section{Introduction}

Let $Z(t)$ be Hardy's function defined by
$$
Z(t)=\zeta(1/2+it) \chi^{-1/2}(1/2+it),
$$
where as usual $\zeta(s)$ is the Riemann zeta-function and $\chi(s)$ is the gamma factor appearing in the functional equation of
$\zeta(s)$:
\begin{equation} \label{FE}
\zeta(s)=\chi(s)\zeta(1-s).
\end{equation}
The explicit form of $\chi(s)$ is
\begin{equation} \label{def-chi}
\chi(s)=2^s\pi^{s-1}\sin\left(\frac{\pi s}{2}\right)\Gamma(1-s)
\end{equation}
and its asymptotic behavour is given by
\begin{equation} \label{chi-asymp}
\chi(\sigma+it)=\left(\frac{|t|}{2\pi}\right)^{1/2-\sigma-it}e^{i(t \pm \frac{\pi}{4})}
 \left(1+O\left(\frac{1}{|t|}\right)\right)
\end{equation}
for $|t| \geq t_0>0$, where $t \pm \frac{\pi}{4}=t+{\rm sgn}(t)\frac{\pi}{4}$. (See Ivi\'c \cite{I0}.)

From \eqref{FE}, it follows that $Z(t)$ is a real-valued even function for real $t$ and $|Z(t)|=|\zeta(1/2+it)|$. 
Therefore the zeros of $\zeta(s)$ on the critical line $\Re s=1/2$ coincide with the real zeros of $Z(t)$. 
Hardy proved that $\int_0^T Z(t)dt \ll T^{7/8}$ and $\int_0^T|Z(t)|dt>\frac12T$, from which he succeeded to show
the infinity of the number of zeros of $\zeta(s)$ on the critical line. (See \cite[p. 51]{Chandra}.)

However, in 2004, Ivi\'c \cite{I2} proved that
$$
\int_0^T Z(t) dt \ll T^{1/4+\varepsilon},
$$
where $\varepsilon$ is an arbitrary small positive number which needs not be the same at each occurrence.
It shows that $Z(t)$ changes sign quite often.
Ivi\'c's result was sharpened by Jutila \cite{J1, J2} and Korolev \cite{K} independently. 
From $Z(t)^2=|\zeta(1/2+it)|^2$, we see that
$$
\int_0^{T} Z^2(t)dt= T \log T +(2\gamma-1-\log 2\pi)T+O(T^{1/3+\varepsilon}).
$$
Since there is a lot of cancellation it is expected that the cubic power moment
has an exponent less than 1. In fact, Ivi\'c showed that
\begin{align*}  
\int_T^{2T}Z^3(t)dt&=2\pi\sqrt{\frac23} \sum_{(\frac{T}{2\pi})^{3/2} \leq n \leq (\frac{T}{\pi})^{3/2}}
\frac{d_3(n)}{n^{1/6}} \cos\left(3\pi n^{2/3}+\frac18\pi\right) \\
& \quad  +O(T^{3/4+\varepsilon})
\end{align*}
and conjectured that
\begin{equation}  \label{conjecture}
 \int_0^{T} Z^3(t)dt \ll T^{3/4+\varepsilon},
\end{equation}
(\cite[Chapter 11]{I4}), but we only know that the left hand side $\ll T (\log T)^{5/2}$ at present.
Here $d_3(n)$ denotes the number of triples $(k_1,k_2, k_3)$ such that $n=k_1k_2k_3, k_j \in \mathbb{Z}, kk_j>0$.

In this paper we shall prove the several mean values of the functions combined with  $Z(t)$ and $\zeta(1/2+it)$.

\begin{theorem}
For large $T>0$, we have
\begin{align*}
\int_0^T Z(t)\zeta\left(\frac12+it\right)dt&=\frac{2\sqrt{2}\pi}{3} e^{\frac{\pi i}{8}}\left(\frac{T}{2\pi}\right)^{3/4}
\left(\frac12 \log \frac{T}{2\pi}+2\gamma-2\log 2 -\frac23 \right) \\
& \quad +O(T^{1/2}\log ^2 T).
\end{align*}
\end{theorem}

\medskip

Ivi\'c's conjecture \eqref{conjecture} would follow from the bound of exponential sum
\begin{equation} \label{es-1}
\sum_{N \leq n \leq 2\sqrt{2}N} \frac{d_3(n)}{n^{1/6}} e^{3\pi i n^{2/3}} \ll N^{1/2+\varepsilon},
\end{equation}
or, as Ivi\'c noted  \cite[(1.6)]{I5}, from
\begin{equation} \label{es-2}
\sum_{N \leq n \leq 2N} d_3(n) e^{3\pi i n^{2/3}} \ll N^{2/3+\varepsilon}.
\end{equation}
It seems that \eqref{es-2} (or \eqref{es-1}) is out of reach of the present method of exponential sums.
However£¬ if we replace $d_3(n)$ by $d(n)$ (the divisor function $d(n)=\sum_{n=d_1d_2}1$), we can prove the
following theorem in the frame of Theorem 1.

\begin{theorem}
Let $A$ be a parameter such that $A \gg N^{-1/4}$.  Then we have
\begin{align*}
&\sum_{N \leq k \leq 2\sqrt{2}N}\frac{d(k)}{k^{1/6}}e^{3\pi i(Ak)^{2/3}} \\
& \quad =\sqrt{3}A^{-4/3} \sum_{A^{4/3}N^{1/3} \leq k \leq \sqrt{2}A^{4/3}N^{1/3}}d(k)k^{1/2}e^{-\pi i (k/A)^2}  \notag \\
& \qquad +O(A^{-1/3}N^{1/2+\varepsilon}) +  O(A^{1/3}N^{1/6} \log N)+O(A^{-1/9}N^{2/9+\varepsilon}) \notag \\[1ex]
& \quad \ll  A^{2/3}N^{1/2} \log N.  \notag
\end{align*}
\end{theorem}

For another kind of mean value of $Z(t)$ and $\zeta(1/2+it)$ we have

\begin{theorem}
For large $T>0$ we have
\begin{align} \label{thm-2}
&\int_0^T Z^2(t) \zeta(1/2+it) dt &= T \left\{\frac12\left(\log \frac{T}{2\pi}\right)^2 + a_1\log\frac{T}{2\pi}+a_2 \right\}
+O(T^{3/4+\varepsilon}),  \notag
\end{align}
where $a_1=3\gamma-1, a_2=3\gamma_1+3\gamma^2-3\gamma+1$, $\gamma_j$ being the coefficients of Laurant expansion of $\zeta(s)$ at
$s=1$ and $\gamma=\gamma_0$ the Euler constant.
\end{theorem}

We note that the integral of the left hand side has an asymptotic form.
It may be interesting to compare with Ivi\'c's conjecture \eqref{conjecture}.

As for another mean value, we shall prove the following

\begin{theorem}
Let  $\alpha$ be a real fixed constant such that $-1/2< \alpha<1/2$. Then we have
$$
\int_T^{2T}Z^3(t) \chi^{\alpha}(1/2+it) dt
 \ll \begin{cases} T^{1-\frac{\alpha}{6}+\varepsilon} & \text{if $ 0 \leq \alpha < 1/2 $}\\[1ex]
                   T^{1+\frac{\alpha}{6}+\varepsilon} & \text{if $ -1/2< \alpha \leq 0$}.
     \end{cases}
$$
\end{theorem}

The cubic moment of Hardy's function corresponds to $\alpha=0$, but unfortunately this gives only $O(T^{1+\varepsilon})$.

%
%
%
%
%


\section{Some lemmas}


%

\begin{lemma} Suppose that $f(x)$ and $\varphi(x)$ are real-valued functions on the interval $[a,b]$
which satisfy the conditions

1) $f^{(4)}(x)$ and $\varphi''(x)$ are continuous.

2) there exist numbers $H,A,U,0<H,A<U,0<b-a \leq U$, such that
$$
A^{-1} \ll f''(x) \ll A^{-1}, \quad f^{(3)} \ll A^{-1}U^{-1}, \quad f^{(4)}(x) \ll A^{-1}U^{-2}
$$
$$
\varphi(x) \ll H, \quad \varphi'(x) \ll HU^{-1}, \quad \varphi''(x) \ll HU^{-2}.
$$

3) $f'(c)=0$ for some $c$, $a \leq c \leq b$.

Then
\begin{align*}
&\int_a^b \varphi(x)\exp(2\pi i f(x)) dx = \frac{1+i}{\sqrt{2}}\frac{\varphi(c)\exp(2\pi if(c)}{\sqrt{f''(c)}}+O(HAU^{-1})\\
& \qquad \quad +O\left(H\min(|f'(a)|^{-1}, \sqrt{A})\right)+O\left(H\min(|f'(b)|^{-1}, \sqrt{A})\right).
\end{align*}
\end{lemma}

This is Lemma 2 of Karatsuba-Voronin \cite[p.71]{KV}. 

\begin{remark}
Here we give an important remark. 
As is noted in Ivi\'c and Zhai \cite{IZ}, the proof actually shows that 
if there is no $c$ which satisfies the condition 3, the term containing $c$ does not appear in the right hand side. 
Moreover if $c=a$ or $c=b$, then the main term is to be halved.
\end{remark}

\begin{lemma}
For $\frac12 \leq  \sigma < 1$ fixed, $1 \ll x, y \ll t^k, s=\sigma+it, xy=(\frac{t}{2\pi})^k, t \geq t_0$ and $k \geq 1$ a fixed integer,
we have
\begin{align*}
\zeta^k(s)&=\sum_{m=1}^{\infty}\rho\left(\frac{m}{x}\right)d_k(m)m^{-s}+\chi^k(s)\sum_{m=1}^{\infty}\rho\left(\frac{m}{y}\right)
d_k(m)m^{s-1} \\
&\quad +O(t^{k(1-\sigma)/3-1})+O(t^{k(1/2-\sigma)-2}y^{\sigma}\log^{k-1}t).  \nonumber
\end{align*}
Here $\chi(s)$ is the function defined by \eqref{def-chi} and $\rho(u) (\geq 0)$ is a smooth function 
such that $\rho(u)+\rho(1/u)=1$ for $u>0$ and $\rho(u)=0$ for $u \geq 2$.  
\end{lemma}

This is Lemma 4 of \cite{IZ}. See also \cite[Theorem 4.16]{I4}.

For the proof of Theorem 4 we need the following lemma.

\begin{lemma}
Let $\alpha, \beta, \gamma$ be fixed real numbers such that $\alpha(\alpha-1)\beta\gamma \neq 0$ and let
$$
S=\sum_{h=H+1}^{2H}\sum_{n=N+1}^{2N} \left|\sum_{M<m \leq 2M}e\left(X\frac{m^{\alpha}h^{\beta}n^{\gamma}}{M^{\alpha}H^{\beta}N^{\gamma}}
\right)\right|^{\ast},
$$
where $\ast$ means that
$$
\left|\sum_{N \leq n \leq N'}z_n\right|^{\ast}=\max_{N \leq N_1 \leq N_2 \leq N'}\left|\sum_{n=N_1}^{N_2}z_n\right|.
$$
Then we have
$$
S \ll (HNM)^{1+\varepsilon} \left\{\left(\frac{X}{HNM^2}\right)^{1/4}+\frac{1}{M^{1/2}}+\frac{1}{X}\right\}.
$$
\end{lemma}

This is Theorem 3 of Robert and Sargos \cite{RS}. Note that $e(x):=e^{2\pi i x}.$


\section{Proofs of Theorem 1 and 2}


\noindent{\it Proof of Theorem 1} \ \
Let $T>0$ be a large number and put
\begin{equation} \label{def-J}
J=\int_T^{2T}Z(t)\zeta\left(\frac12+it\right)dt.
\end{equation}
By the definition of $Z(t)$ and applying Lemma 2 we have
\begin{align*}
Z(t)\zeta\left(\frac12+it\right)&=\zeta^2\left(\frac12+it\right)\chi^{-1/2}\left(\frac12+it\right) \\
&=\left(\sum_{k=1}^{\infty}\rho\left(\frac{k}{x}\right)\frac{d(k)}{k^{1/2+it}}+\chi^2\left(\frac12+it\right)
\sum_{k=1}^{\infty}\rho\left(\frac{k}{y}\right)\frac{d(k)}{k^{1/2-it}} \right. \\
& \qquad + O\left(t^{-2/3}\right)+O\left(t^{-2}y^{1/2} \log t \right) \biggl) \chi^{-1/2}\left(\frac12+it\right),
\end{align*}
where $xy=(t/2\pi)^2$.
Substituting this expression to \eqref{def-J}, we have
\begin{equation} \label{J-bunkai}
J=J_1+J_2+O(T^{1/3}),
\end{equation}
where
\begin{align}
J_1&=\sum_{k=1}^{\infty}\frac{d(k)}{k^{1/2}}\int_T^{2T}\rho\left(\frac{k}{x}\right)k^{-it}\chi^{-1/2}\left(\frac12+it\right)dt, \label{J1form}  \\
\intertext{and}
J_2&=\sum_{k=1}^{\infty}\frac{d(k)}{k^{1/2}}\int_T^{2T}\rho\left(\frac{k}{y}\right)k^{it}\chi^{3/2}\left(\frac12+it\right)dt.
\label{J2form}
\end{align}

We take
$$
x=2\left(\frac{t}{2\pi}\right), \quad y=\frac12 \left(\frac{t}{2\pi}\right),
$$
and put $K=\frac{T}{\pi}.$  Then the ranges of $k$ in the sums in \eqref{J1form} and \eqref{J2form} 
are in fact $k \leq 4K $ and $k \leq K$ respectively. 

We first consider $J_1$.  By \eqref{chi-asymp}, we find that
$$
k^{-it}\chi^{-1/2}\left(\frac12+it\right)=e^{-\frac{\pi i}{8}}e^{\frac{i}{2}(\tlg-t\log k^2)}+O(1/t),
$$
therefore we have
$$
J_1=e^{-\frac{\pi i}{8}} \sum_{k \leq 4K}\frac{d(k)}{k^{1/2}}\int_T^{2T}
\rho\left(\frac{k}{x}\right)e^{\frac{i}{2}(\tlg-t\log k^2)}dt +O(T^{1/2}\log T).
$$
We evaluate the integral by applying Lemma 1
with $\varphi(t)=\rho\left(k\left(\frac{\pi}{t}\right)\right), f(t)=\frac{1}{4\pi}(\tlg-t\log k^2)$.
Note that $\varphi(t)$ satisfies the conditions of Lemma 1 with $H=1, U=T$.
Since $f'(t_0)=0$ if and only if $t_0=2\pi k^2$, the main term of the integral appears for $k$ such that
\begin{equation} \label{k-hanni}
\left(\frac{T}{2\pi}\right)^{1/2} \leq k \leq \left(\frac{T}{\pi}\right)^{1/2}.
\end{equation}
Hence we get
\begin{align*}
&\int_T^{2T}\rho\left(\frac{k}{x}\right)e^{\frac{i}{2}(\tlg-t\log k^2)}dt \\
&=M(k)+O\left(1+\min\Bigl(\sqrt{T}, \frac{1}{|\log(\frac{T}{2\pi k^2})|}\Bigr)+
                \min\Bigl(\sqrt{T}, \frac{1}{|\log(\frac{T}{\pi k^2})|}\Bigr) \right),  \notag
\end{align*}
where
$$ M(k)= e^{\frac{\pi i}{4}}\rho\left(\frac{1}{2k}\right) 2\sqrt{2} \pi k e^{-\pi i k^2}=2\sqrt{2} \pi e^{\frac{\pi i }{4}}k (-1)^k$$
for $k$ satisfying the condition \eqref{k-hanni} and 0 otherwise.  This yields that
\begin{align*}
J_1&=2\sqrt{2}\pi e^{\frac{\pi i}{8}}\sum_{(\frac{T}{2\pi})^{1/2} \leq k \leq (\frac{T}{\pi})^{1/2}} \hspace{-8mm}' \qquad
            (-1)^kd(k)k^{1/2} \\
   & \quad  +\sum_{k \leq 4K}\frac{d(k)}{k^{1/2}}\ O\left(1+ \min\Bigl(\sqrt{T}, \frac{1}{|\log(\frac{T}{2\pi k^2})|}\Bigr)+
                \min\Bigl(\sqrt{T}, \frac{1}{|\log(\frac{T}{\pi k^2})|}\Bigr)\right) \\
   & \quad +O(T^{1/2}\log T) \\
   &=:R_0+R_1+R_2+R_3 +O(T^{1/2}\log T),
\end{align*}
where $\sum'$ means that the terms for $k=(T/2\pi)^{1/2}$ and $k=(T/\pi)^{1/2}$ are to be halved if they are 
integers.  It is clear that $R_1 \ll T^{1/2} \log T. $
To estimate $R_2$, we divide the sum into four parts: 
\begin{align*}
\sum_{k \leq 4K}&=\sum_{1 \leq k < \frac12(\frac{T}{2\pi})^{1/2}}
+\sum_{\frac12(\frac{T}{2\pi})^{1/2} \leq k < (\frac{T}{2\pi})^{1/2}}
+\sum_{(\frac{T}{2\pi})^{1/2} \leq k \leq  2(\frac{T}{2\pi})^{1/2}}
+ \sum_{2(\frac{T}{2\pi})^{1/2} < k \leq 4K} \\
&=:S_1+S_2+S_3+S_4.
\end{align*}
For $S_1$, since $\min(\sqrt{T}, \frac{1}{|\log(\frac{T}{2\pi k^2})|}) \ll \frac{1}{\log 4}$, we have $S_1 \ll T^{1/4}\log T$.
For $S_4$, we have the same upper bound for $\min(\sqrt{T}, \frac{1}{|\log(\frac{T}{2\pi k^2})|})$,  
hence we have $S_4 \ll T^{1/2} \log T$. 
Now we consider $S_2$. 
We write $k=[(\frac{T}{2\pi})^{1/2}]-j$ for $k$ in this range and set $S_2=S_{2,1}+S_{2,2}$, where
$S_{2,1}$ is the sum for $j=0,1,2$ and $S_{2,2}$ is the sum for $j \geq 3$.  
For $S_{2,1}$ we adopt $\min(\sqrt{T}, \frac{1}{|\log(\frac{T}{2\pi k^2})|})= \sqrt{T}$, 
hence $S_{2,1} \ll T^{1/4+\varepsilon}$. For $S_{2,1}$, we have 
$$ 
\log \frac{(\frac{T}{2\pi})^{1/2}}{k}=\left|\log\frac{[(\frac{T}{2\pi})^{1/2}]-j}{(\frac{T}{2\pi})^{1/2}}\right| 
\asymp \frac{j}{(\frac{T}{2\pi})^{1/2}},
$$ 
from which we get
$$
S_{2,2} \ll \sum_{j} \frac{d(k)}{k^{1/2}}\frac{(\frac{T}{2\pi})^{1/2}}{j} \ll T^{1/4+\varepsilon}
$$
Therefore $S_2 \ll T^{1/4+\varepsilon}$. It is the same for $S_3$.
Taken together we have $ R_2 \ll T^{1/2}\log T.$  Similarly we have $ R_3 \ll T^{1/2}\log T.$  

%

 As  a result, we get
\begin{equation}  \label{Th1-J1-hyouka}
J_1 = 2\sqrt{2}\pi e^{\frac{\pi i}{8}}\sum_{(\frac{T}{2\pi})^{1/2} \leq k \leq (\frac{T}{\pi})^{1/2}}
       \hspace{-8mm}' \quad (-1)^kd(k)k^{1/2}+O(T^{1/2}\log  T).
\end{equation}

\medskip

Next we consider $J_2$.  Similarly to the case of $J_1$, we have by \eqref{chi-asymp},
\begin{align}  \label{J2}
J_2& =e^{\frac{3\pi i}{8}}\sum_{k \leq K} \frac{d(k)}{k^{1/2}}
\int_T^{2T}\rho\left(\frac{k}{y}\right)e^{-\frac32 i (\tlg-t\log k^{2/3})}dt  \\
& \quad + O(T^{1/2}\log T). \notag
\end{align}
We apply Lemma 1 to the above integral with $\varphi(t)=\rho(2k(2\pi/t))$ and $f(t)=-\frac{3}{4\pi}(\tlg-t\log k^{2/3})$.
In this case $f'(t_0)=0$ if and only if $t_0=2\pi k^{2/3}$ and $t_0$ is contained in the interval $[T,2T]$ if and only if
$(\frac{T}{2\pi})^{3/2} \leq k \leq (\frac{T}{\pi})^{3/2}$. Since the range of the sum over $k$ is $1 \leq k \leq K$,
there are no such $k$, that is, the integral in \eqref{J2} does not have the main term.
Considering the error term by Lemma 1 we find that
\begin{align*}
J_2 & \ll \sum_{k \leq K}\frac{d(k)}{k^{1/2}}\left(1+\min\Bigl((\sqrt{T},\frac{1}{|\log\frac{T}{2\pi k^{2/3}}|}\Bigr)
              +\min\Bigl(\sqrt{T}, \frac{1}{|\log\frac{T}{\pi k^{2/3}}|}\Bigr)\right) \\
& =: R_1'+R_2'+R_3'.
\end{align*}
We have clearly $R_1' \ll T^{1/2}\log T$. For $R_2'$ and $R_3'$ we note that $|\log \frac{T}{k^{2/3}}| \gg 1 $
since $k \leq K$, which implies that $R_2', R_3' \ll T^{1/2}\log T$. Hence
\begin{equation}  \label{Th1-J2-hyouka}
J_2 \ll T^{1/2} \log T.
\end{equation}

From \eqref{J-bunkai}, \eqref{Th1-J1-hyouka} and \eqref{Th1-J2-hyouka}, we get
$$
J = 2\sqrt{2}\pi e^{\frac{\pi i}{8}}\sum_{(\frac{T}{2\pi})^{1/2} \leq k \leq (\frac{T}{\pi})^{1/2}}
    \hspace{-8mm}' \quad (-1)^kd(k)k^{1/2}+O(T^{1/2}\log  T).
$$

Now dividing the interval $[0,T]$ as $ \cup_j [T/2^j, T/2^{j-1}]$ and summing the above evaluations we see that
\begin{align} \label{Th1-0}
\int_0^T Z(t)\zeta\left(\frac12+it\right)dt
&=2\sqrt{2}\pi e^{\frac{\pi i}{8}}\sum_{k \leq (\frac{T}{2\pi})^{1/2}}(-1)^kd(k)k^{1/2} \\
& \quad +O(T^{1/2}\log^2 T).  \nonumber 
\end{align}
It is known that  for $x \gg 1$
$$
\sum_{k \leq x}(-1)^k d(k)=\frac{x}{2}(\log x+2\gamma-1-2\log 2)+O(x^{1/3+\varepsilon}),
$$
(see e.g. Ivic \cite{I3}). By partial summation we get
$$
\sum_{k \leq x}(-1)^k d(k)k^{1/2}=\frac13 x^{3/2} \left(\log x+2\gamma-2\log 2-\frac23\right)+O(x^{5/6+\varepsilon}).
$$
Substituting this form in \eqref{Th1-0} we get
\begin{align*}
\int_0^T Z(t)\zeta\left(\frac12+it\right)dt&=\frac{2\sqrt{2}\pi}{3} e^{\frac{\pi i}{8}}\left(\frac{T}{2\pi}\right)^{3/4}
\left(\frac12 \log \frac{T}{2\pi}+2\gamma-2\log 2 -\frac23 \right) \\
& \quad +O(T^{1/2}\log ^2 T).
\end{align*}
This proves the assertion of Theorem 1.

\bigskip

\noindent{\it Proof of Theorem 2.} \ \
Let $A$ be a parameter such that $T^{-1/2} \ll A \ll T^{3/2}$.  
We shall consider the integral
\begin{equation*}
J_A=\int_T^{2T}Z(t)\zeta\left(\frac12+it\right)A^{it}dt
\end{equation*}
by the same way as in the proof of Theorem 1.  Applying Lemma 2 we get
\begin{equation}  \label{JA-bunkai}
J_A= J_{A,1}+J_{A,2}+O(T^{1/3}),
\end{equation}
where we put
\begin{align}
J_{A,1}=\int_T^{2T} \chi^{-1/2}\left(\frac12+it\right)\sum_{k=1}^{\infty}
                     \rho\left(\frac{k}{x}\right)\frac{d(k)}{k^{1/2+it}}A^{it}dt \label{def-JA1} \\
\intertext{and}
J_{A,2}=\int_T^{2T} \chi^{3/2}\left(\frac12+it\right)\sum_{k=1}^{\infty}
                     \rho\left(\frac{k}{y}\right)\frac{d(k)}{k^{1/2-it}}A^{it}dt, \label{def-JA2}
\end{align}
where $xy=(\frac{t}{2\pi})^2$. We shall evaluate both $J_{A,1}$ and $J_{A,2}$ by taking two different choices of $x$ and $y$.
Hereafter we put 
$$
K_0=\left(\frac{T}{\pi}\right)^{1/2}.
$$

\medskip

\noindent{\bf The case $x=8A(\frac{t}{2\pi})^{1/2}$ and $y=\frac{1}{8A}(\frac{t}{2\pi})^{3/2}$.}
The ranges of $J_{A,1}$ and $J_{A,2}$ are  at most $k \leq 16AK_0$ and $k \leq \frac{1}{4A}K_0^3$, respectively.
By \eqref{chi-asymp} and the trivial estimate for the error term we get
\begin{align}  \label{JA-1}
J_{A,1}&=e^{-\frac{\pi i}{8}}\sum_{k \leq 16AK_0} \frac{d(k)}{k^{1/2}}\int_T^{2T} \rho\left(\frac{k}{x}\right)
        e^{\frac{i}{2}(\tlg-t\log(\frac{k}{A})^2)} dt \\
       & \quad + O\left(A^{1/2}T^{1/4+\varepsilon}\right),  \notag
\end{align}
We shall evaluate the integral by Lemma 1. Let $f(t)=\frac{1}{4\pi}(\tlg-t\log(\frac{k}{A})^2)$. Then $f'(t_0)=0$ if and
only if $t_0=2\pi(\frac{k}{A})^2$ and $T \leq t_0 \leq 2T$ if and only if
\begin{equation} \label{k-range-JA1}
A\left(\frac{T}{2\pi}\right)^{1/2} \leq k \leq A\left(\frac{T}{\pi}\right)^{1/2}.
\end{equation}
We see that all $k$ satisfing \eqref{k-range-JA1} are contained
in the range $k \leq 16AK_0$. Therefore the integral in \eqref{JA-1} has a main term which is given by
$$
M_A(k)= e^{\frac{\pi i}{4}}\rho\left(\frac18\right)2\sqrt{2}\pi \frac{k}{A} e^{-\pi i (k/A)^2}
$$
for $A(\frac{T}{2\pi})^{1/2} \leq k \leq A(\frac{T}{\pi})^{1/2}$ and $ M_A(k)=0 $ otherwise.
We note that $\rho(1/8)=1$ in the above formula.
It follows from Lemma 1 and \eqref{JA-1} that

\begin{align*}
&J_{A,1}=e^{-\frac{\pi i}{8}}\sum_{A(\frac{T}{2\pi})^{1/2} \leq k \leq A(\frac{T}{\pi})^{1/2}}\frac{d(k)}{k^{1/2}}M_A(k) \\
 & + \sum_{k \leq 4AK_0} \frac{d(k)}{k^{1/2}} \ O\left( 
          1+ \min\biggl(\sqrt{T}, \frac{1}{|\log(\frac{(T/2\pi)^{1/2}}{k/A})|}\biggr)
           + \min\biggl(\sqrt{T}, \frac{1}{|\log(\frac{(T/\pi)^{1/2}}{k/A})|}\biggr)\right) \\[1ex]
 & +O(A^{1/2}T^{1/4+\varepsilon}). 
\end{align*}
%
%
%
%
%
%
%
%
Similarly to the proof of Theorem 1, we see that the contributions from the $O$-terms are bounded by
$O(A^{1/2}T^{1/4+\varepsilon}+A^{-1/2}T^{1/4+\varepsilon})$.  
%
Hence we get
\begin{align}  \label{JA1}
J_{A,1}& =e^{\frac{\pi i}{8}}\frac{2\sqrt{2}\pi}{A}\sum_{A(\frac{T}{2\pi})^{\frac12}\leq k \leq A(\frac{T}{\pi})^{\frac12}}
          d(k)k^{1/2}e^{-\pi i (k/A)^2} \\[1ex]
       & \quad {} +O(A^{1/2}T^{1/4+\varepsilon})+O(A^{-1/2}T^{1/4+\varepsilon}).  \nonumber 
\end{align}

Next we consider $J_{A,2}$. Similarly to $J_{A,1}$ we have
\begin{align*}
J_{A,2}&=e^{\frac{3\pi i}{8}}\sum_{k \leq \frac{1}{4A}K_0^3}\frac{d(k)}{k^{1/2}}\int_T^{2T} \rho\left(\frac{k}{y}\right)
        e^{-\frac{3}{2}i(\tlg-t\log(Ak)^{2/3})}dt \\
& \quad +O(A^{-1/2}T^{3/4+\varepsilon}).
\end{align*}
If we put $f(t)=-\frac{3}{4\pi}(\tlg-t\log(Ak)^{2/3})$ this time, $f'(t_0)=0$ if and only if
$t_0=2\pi(Ak)^{2/3}$ and so $T \leq t_0 \leq 2T$ if and only if
\begin{equation} \label{k-range-JA2}
\frac{1}{A}\left(\frac{T}{2\pi}\right)^{3/2} \leq k \leq \frac{1}{A}\left(\frac{T}{\pi}\right)^{3/2}.
\end{equation}
Since $K$ runs over $1 \leq k \leq \frac{1}{4A}K_0^3$ there is no main term in the integral of $J_{A,2}$. 
Hence by Lemma 1, 
we get similarly that
\begin{align}  \label{JA2}
J_{A,2} & \ll \sum_{k \leq \frac{1}{4A}K_0^3}\frac{d(k)}{k^{1/2}} 
          \left(1+\min\biggl(\sqrt{T}, \frac{1}{|\log(\frac{(T/2\pi)^{3/2}}{Ak})|}\biggr) \right.  \\
        & \hspace{3.3cm} \left.  +\min(\sqrt{T}, \frac{1}{|\log(\frac{(T/\pi)^{3/2}}{Ak})|}\biggr)\right) \nonumber \\
        & \ll A^{-1/2}T^{3/4+\varepsilon}+A^{1/2}T^{-1/4+\varepsilon}.   \nonumber
\end{align}

From \eqref{JA-bunkai}, \eqref{JA1} and \eqref{JA2}, we obtain
\begin{align}  \label{JA-firstresult}
J_{A}& =e^{\frac{\pi i}{8}}\frac{2\sqrt{2}\pi}{A}\sum_{A(\frac{T}{2\pi})^{\frac12}\leq k \leq A(\frac{T}{\pi})^{\frac12}}
          d(k)k^{1/2}e^{-\pi i (k/A)^2} \\
     & \quad {}+O(A^{1/2}T^{1/4+\varepsilon})+O(A^{-1/2}T^{3/4+\varepsilon})+O(T^{1/3}). \notag  
\end{align}

\medskip

\noindent{\bf The case $x=\frac{A}{4}(\frac{t}{2\pi})^{1/2}$ and $y=\frac{4}{A}(\frac{t}{2\pi})^{3/2}$.} 
In this choice of $x$ and $y$, the sums in \eqref{def-JA1} and \eqref{def-JA2} are actually 
over $k \leq \frac12 A K_0$ and $k \leq \frac{8}{A}{K_0}^3$ respectively.
Thus
\begin{align*}  
J_{A,1}&=e^{-\frac{\pi i}{8}}\sum_{k \leq \frac A2 K_0} \frac{d(k)}{k^{1/2}}\int_T^{2T}
       \rho\left(\frac{k}{x}\right)e^{\frac{i}{2}(\tlg-t\log(\frac{k}{A})^2)}dt \\
       & \quad +O\left(A^{1/2}T^{1/4+\varepsilon}\right) \nonumber \\
\intertext{and} 
J_{A,2}&=e^{\frac{3\pi i}{8}}\sum_{k \leq \frac 8A K_0^3} \frac{d(k)}{k^{1/2}}\int_T^{2T}
      \rho\left(\frac{k}{y}\right)e^{-\frac{3}{2}i(\tlg-t\log(Ak)^{2/3})}dt  \\
       & \quad + O(A^{-1/2}T^{3/4+\varepsilon}). \nonumber
\end{align*}

As for $J_{A,1}$, the integral has a main term if and only if $k$ satisfies \eqref{k-range-JA1}.
Since $k$ runs over $1 \leq k \leq \frac A2 K_0$, there are no such $k$. 
The contribution from the error term of the integral is the same as in the previous case since the range of the sum
has the same order, hence we get
\begin{equation}  \label{JA1-second}
J_{A,1} \ll A^{1/2}T^{1/4+\varepsilon}+A^{-1/2}T^{1/4+\varepsilon}.
\end{equation}

%

On the other hand, the integral of $J_{A,2}$ has a main term if and only if $k$ satisfies \eqref{k-range-JA2},
and in fact all $k$ are in the range $k \leq \frac{8}{A}K_0^3$. Hence by Lemma 1,
 $J_{A,2}$ has the following form:
\begin{align*} 
J_{A,2}&=e^{\frac{3\pi i}{8}}\sum_{\frac{1}{A}(\frac{T}{2\pi})^{3/2} \leq k \leq \frac1A(\frac{T}{\pi})^{3/2}}
         \frac{d(k)}{k^{1/2}}\widetilde{M}_A(k) \\
       & \quad + \sum_{k \leq \frac8A K_0^3} \frac{d(k)}{k^{1/2}} 
               \ O\left(1+\min\biggl(\sqrt{T}, \frac{1}{|\log(\frac{(T/2\pi)^{3/2}}{Ak})|}\biggr)  \right. \nonumber \\
       & \hspace{3.8cm} \left.+\min(\sqrt{T}, \frac{1}{|\log(\frac{(T/\pi)^{3/2}}{Ak})|}\biggr)  \right)  \nonumber \\
& \quad +O(A^{-1/2}T^{3/4+\varepsilon}),  \nonumber
\end{align*}
where
$$
\widetilde{M}_A(k)=e^{-\frac{\pi i}{4}}\rho\left(\frac14\right)\frac{2\sqrt{2}\pi}{\sqrt{3}}(Ak)^{1/3}e^{3\pi i (Ak)^{2/3}}
$$
for $\frac{1}{A}(\frac{T}{2\pi})^{3/2} \leq k \leq \frac1A(\frac{T}{\pi})^{3/2}$ and 0 otherwise.
%
%
We see that the contribution from the $O$-term is the same as the previous case, therefore
\begin{align}  \label{JA2-second}
J_{A,2}&=e^{\frac{\pi i}{8}}\frac{2\sqrt{2}\pi}{\sqrt{3}}A^{1/3}
\sum_{\frac{1}{A}(\frac{T}{2\pi})^{3/2} \leq k \leq \frac1A(\frac{T}{\pi})^{3/2}}\frac{d(k)}{k^{1/6}}e^{3\pi i(Ak)^{2/3}} \\
&\quad + O(A^{-1/2}T^{3/4+\varepsilon})+O(A^{1/2}T^{-1/4+\varepsilon}).  \nonumber
\end{align}

From \eqref{JA1-second} and \eqref{JA2-second} we obtain that
\begin{align}  \label{JA-secondresult}
J_A&=e^{\frac{\pi i}{8}}\frac{2\sqrt{2}\pi}{\sqrt{3}}A^{1/3}
\sum_{\frac{1}{A}(\frac{T}{2\pi})^{3/2} \leq k \leq \frac1A(\frac{T}{\pi})^{3/2}}\frac{d(k)}{k^{1/6}}e^{3\pi i(Ak)^{2/3}} \\
&\quad + O(A^{-1/2}T^{3/4+\varepsilon})+O(A^{1/2}T^{1/4+\varepsilon})+O(T^{1/3}).  \nonumber
\end{align}

Now we have two expressions of $J_A$: \eqref{JA-firstresult} and \eqref{JA-secondresult}.  Comparing these expressions we obtain
\begin{align}  \label{iikae-1}
&\sum_{\frac{1}{A}(\frac{T}{2\pi})^{3/2} \leq k \leq \frac1A(\frac{T}{\pi})^{3/2}}\frac{d(k)}{k^{1/6}}e^{3\pi i(Ak)^{2/3}} \\[1ex]
& \hspace{1cm} =\sqrt{3}A^{-4/3} \sum_{A(\frac{T}{2\pi})^{\frac12}\leq k \leq A(\frac{T}{\pi})^{\frac12}}d(k)k^{1/2}e^{-\pi i (k/A)^2} \notag \\
& \hspace{1cm} \quad +O(A^{-5/6}T^{3/4+\varepsilon}) +O(A^{1/6}T^{1/4+\varepsilon})+O(A^{-1/3}T^{1/3+\varepsilon}) \notag \\[1ex]
& \hspace{1cm} \ll  A^{1/6}T^{3/4} \log T,  \notag
\end{align}
where the last inequality is obtained by the trivial estimate.
In \eqref{iikae-1}, we take $T=2\pi (AN)^{2/3}$. Then \eqref{iikae-1} is transformed to
\begin{align*} 
&\sum_{N \leq k \leq 2\sqrt{2}N}\frac{d(k)}{k^{1/6}}e^{3\pi i(Ak)^{2/3}} \\
& \quad =\sqrt{3}A^{-4/3} \sum_{A^{4/3}N^{1/3} \leq k \leq \sqrt{2}A^{4/3}N^{1/3}}d(k)k^{1/2}e^{-\pi i (k/A)^2}  \notag \\
& \qquad +O(A^{-1/3}N^{1/2+\varepsilon}) +  O(A^{1/3}N^{1/6+\varepsilon})+O(A^{-1/9}N^{2/9+\varepsilon}) \notag \\[1ex]
& \quad \ll  A^{2/3}N^{1/2} \log N  \notag
\end{align*}
for $A \gg N^{-1/4}$.  This proves the assertion of Theorem 2.


\section{Proof of Theorem 3}


The method is similar to the previous cases, but we shall write the necessary points for the sake of
completeness. Let $T$ be a large number.  We put  
\begin{equation*}
I=\int_T^{2T} Z^2(t)\zeta\left(\frac12+it\right)dt.
\end{equation*}
By the definition of Hardy's function and Lemma 2, we have
\begin{align} \label{AFE-3}
&Z^2(t)\zeta\left(\frac12+it\right)=\zeta^3\left(\frac12+it\right)\chi^{-1}\left(\frac12+it\right) \\
&=\chi^{-1}\left(\frac12+it\right)\sum_{k=1}^{\infty}\rho\left(\frac{k}{x}\right)\frac{d_3(k)}{k^{1/2+it}}
+ \chi^2\left(\frac12+it\right)\sum_{k=1}^{\infty}\rho\left(\frac{k}{y}\right)\frac{d_3(k)}{k^{1/2-it}} \notag \\
& \quad + O\left(t^{-1/2}\right)+O\left(t^{-2}y^{1/2}\log^ 2 t\right), \notag
\end{align}
where $xy=(\frac{t}{2\pi})^3$.

We take $x=2(\frac{t}{2\pi})^{3/2}$ and $y=\frac12(\frac{t}{2\pi})^{3/2}$ in \eqref{AFE-3} and 
put $K_3=(T/\pi)^{3/2}$. Then the ranges of $k$ in the above two sums are at most $k \leq 4K_3$ and $k \leq K_3$, respectively.
Hence
\begin{align}  \label{thm3-I}
I&=\sum_{k \leq 4K_3}\frac{d_3(k)}{k^{1/2}}\int_T^{2T}\rho\left(\frac{k}{x}\right)k^{-it}\chi^{-1}\left(\frac12+it\right)dt \\
& \quad +\sum_{k \leq K_3}\frac{d_3(k)}{k^{1/2}}\int_T^{2T}\rho\left(\frac{k}{y}\right)k^{it}\chi^{2}\left(\frac12+it\right)dt
 + O(T^{1/2}) \notag \\
&=: I_1+I_2+O(T^{1/2}).  \notag
\end{align}

As for $I_2$, using \eqref{chi-asymp}, we get
\begin{align*}
I_2
&=e^{\frac{\pi i}{2}}\sum_{k \leq K_3}\frac{d_3(k)}{k^{1/2}}
      \int_T^{2T} \rho\left(\frac{k}{y}\right) e^{-2i(\tlg-t\log\sqrt{k})} dt+ O(T^{3/4}\log^2 T).
\end{align*}
As previously, we apply Lemma 1 to the above integral with $\varphi(t)=\rho\left(2k \left(\frac{2\pi}{t}\right)^{3/2}\right)$ and
$f(t)=-\frac{1}{\pi}(\tlg-t\log\sqrt{k}) $.
We see that $f'(t_0)=0$ if and only if $t_0=2\pi \sqrt{k}$,
and this $t_0$ is contained in the interval $[T, 2T]$ if and only if
\begin{equation} \label{I2-k-hanni}
\left(\frac{T}{2\pi}\right)^2 \leq k \leq \left(\frac{T}{\pi}\right)^2.
\end{equation}
Since $k$ runs over the range $1 \leq k \leq K_3$, there is no $k$ which satisfies \eqref{I2-k-hanni}, hence
the main term does not appear in this integral. On the other hand, the error term of this integral is given by
$ 1+ \min\left(\sqrt{T}, \frac{1}{|\log\frac{T}{2\pi\sqrt{k}}|}\right)
+\min\left(\sqrt{T}, \frac{1}{|\log\frac{T}{\pi\sqrt{k}}|}\right) \ll 1$,
hence we get
\begin{align}  \label{I2}
I_2 \ll \sum_{k \leq K_3} \frac{d_3(k)}{k^{1/2}} \ll T^{3/4} \log^2 T.
\end{align}

Next we treat $I_1$. By \eqref{chi-asymp} again, we have
\begin{align*}
I_1 
   &=e^{-\frac{\pi i}{4}}\sum_{k \leq 4K_3}\frac{d_3(k)}{k^{1/2}}
        \int_T^{2T} \rho\left(\frac{k}{x}\right) e^{i(\tlg-t\log k)} dt+O(T^{3/4}\log^2 T).
\end{align*}
In this case $\varphi(t)=\rho(k(\frac{2\pi}{t})^{3/2}/2)$ and $f(t)=\frac{1}{2\pi}(\tlg-t\log k)$.
We see that $f'(t_0)=0$  if and only if $t_0=2\pi k$ and this $t_0$ is contained in $[T,2T]$ if and only if
\begin{equation*}  
\frac{T}{2\pi} \leq k \leq \frac{T}{\pi}.
\end{equation*}
Hence we have
\begin{align*}
&\int_T^{2T}\rho\left(\frac{k}{x}\right)e^{i(\tlg-t\log k)}dt \\
&=M(k)+O\left(1+\min\left(\sqrt{T}, \frac{1}{|\log\frac{T}{2\pi k}|}\right)+\min\left(\sqrt{T}, \frac{1}{|\log\frac{T}{\pi k}|}\right)\right),
\end{align*}
where $M(k)$ is the main term given by
\begin{align*}
M(k)&=e^{\frac{\pi i}{4}}\rho\left(\frac{1}{2\sqrt{k}}\right) (2\pi t_0)^{1/2} e^{-2\pi i k} =2\pi e^{\frac{\pi i }{4}}k^{1/2}
\end{align*}
for $k$ such that $\frac{T}{2\pi} \leq k \leq \frac{T}{\pi}$ and $0$ otherwise. 
Therefore we get
\begin{align}  \label{I1}
I_1&=2\pi \sum_{\frac{T}{2\pi} \leq k \leq \frac{T}{\pi}}{\hspace{-3mm} '} \quad d_3(k) \\
   &\quad + \sum_{k \leq 4K_3}\frac{d_3(k)}{k^{1/2}}\left(1+\min\left(\sqrt{T}, \frac{1}{|\log\frac{T}{2\pi k}|}\right)
           +\min\left(\sqrt{T}, \frac{1}{|\log\frac{T}{\pi k}|}\right)\right) \nonumber \\
   & =  2\pi \sum_{\frac{T}{2\pi} \leq k \leq \frac{T}{\pi}}{\hspace{-3mm} '} \quad d_3(k)+O(T^{3/4} \log^2T).  \notag
\end{align}
Here we can get the last $O$-term by the same way as previously.
%
Combining \eqref{thm3-I}, \eqref{I2} and \eqref{I1}, we obtain
\begin{align*}
I=2\pi \sum_{\frac{T}{2\pi} \leq k \leq \frac{T}{\pi}} {\hspace{-3mm} '} \quad d_3(k)+O(T^{3/4}\log^2T).
\end{align*}

Now dividing the interval $[0,T]$ as $ \cup_j [2^{j}T, 2^{j-1}T]$ we obtain that
$$
\int_0^T Z^2(t) \zeta\left(\frac12+it\right)dt= 2\pi \sum_{k \leq \frac{T}{2\pi}}d_3(k)+O(T^{3/4}\log^3 T).
$$
Theorem 3 follows from the well-known formula:
$$
\sum_{n \leq x}d_3(n)=x\left(\frac12 \log^2 x+(3\gamma-1)\log x+3\gamma_1+3\gamma^2-3\gamma+1\right)+O(x^{1/2}),
$$
where $\gamma_j$ is the coefficients of Laurant expansion of $\zeta(s)$ at $s=1$.


\section{Proof of Theorem 4}


Let
$$
I=\int_T^{2T}Z^3\left(\frac12+it\right)\chi^{\alpha}\left(\frac12+it\right)dt,
$$
where $\alpha$ is a fixed constant such that $-1/2<\alpha<1/2$. By the definition of $Z(t)$ and Lemma 2 we have
\begin{align*}
I=I_1+I_2+O(T^{1/2}),
\end{align*}
where
\begin{align}
I_1&=\sum_{k=1}^{\infty}\frac{d_3(k)}{k^{1/2}}\int_T^{2T}\rho\left(\frac{k}{x}\right)k^{-it}
      \chi^{\alpha-\frac32}\left(\frac12+it\right)dt \label{Th4-I1} \\
\intertext{and}
I_2&=\sum_{k=1}^{\infty}\frac{d_3(k)}{k^{1/2}}\int_T^{2T}\rho\left(\frac{k}{y}\right)k^{it}
      \chi^{\alpha+\frac32}\left(\frac12+it\right)dt,  \label{Th4-I2}
\end{align}
where $xy=(\frac{t}{2\pi})^3$. The evaluations of these integrals are the same as before, so we only sketch the outline of
these evaluations.

Assume that $0 \leq \alpha < \frac12$. We take  $x=2(\frac{t}{2\pi})^{1/2}$ and $y=\frac12(\frac{t}{2\pi})^{1/2}$ 
and put $K_4=(\frac{T}{\pi})^{3/2}$.
Then $k$ in the summations in (\ref{Th4-I1}) and (\ref{Th4-I2}) at most run over $1 \leq k \leq 4K_4$
 and $1 \leq k \leq K_4$ respectively.

We shall treat $I_1$ first. By \eqref{chi-asymp}, we see that the integral in (\ref{Th4-I1}) becomes
$$
e^{\frac{\pi i}{4}(\alpha-\frac32)}\int_T^{2T}\rho\left(\frac{k}{x}\right)
e^{(\frac32-\alpha)i(\tlg-t\log k^{\frac{1}{3/2-\alpha}})}dt +O(1).
$$
The main term of the integral above appears only when
\begin{equation*} 
\left(\frac{T}{2\pi}\right)^{\frac32-\alpha} \leq k \leq \left(\frac{T}{\pi}\right)^{\frac32-\alpha},
\end{equation*}
in which case it is given by
$$
M_{\alpha}(k)=e^{\frac{\pi i}{4}}\rho\left(k^{\frac{-2\alpha}{3-2\alpha}}/2\right)\frac{2\pi}{\sqrt{3/2-\alpha}}\ k^{\frac{1}{3-2\alpha}}
\ e^{-(\frac32-\alpha)i k^{\frac{1}{3/2-\alpha}}}.
$$
Computing the error term by Lemma 1, we get
\begin{align}  \label{Th4-I1-main-sum}
I_1&=e^{\frac{\pi i}{4}(\alpha-\frac12)}\frac{2\pi}{\sqrt{3/2-\alpha}}
    \sum_{(\frac{T}{2\pi})^{3/2-\alpha} \leq k \leq (\frac{T}{\pi})^{3/2-\alpha}} \frac{d_3(k)}{k^{1/2}}k^{\frac{1}{3-2\alpha}}e^{-(\frac32-\alpha)ik^{1/(3/2-\alpha)}} \\
   & \quad + \sum_{k \leq K_4}\frac{d_3(k)}{k^{1/2}}
     \ O \left(1+\min\biggl(\sqrt{T}, \frac{1}{|\log\frac{(T/2\pi)^{3/2-\alpha}}{k}|}\biggr) \right. \nonumber \\
   & \hspace{5cm} \left. +\min\biggl(\sqrt{T}, \frac{1}{|\log\frac{(T/\pi)^{3/2-\alpha}}{k}|}\biggr) \right). \notag
\end{align}
Just in the same way in the previous cases, we can see easily that the above $O$-term is estimated as $O(T^{3/4}\log^2 T)$.

On the other hand, for $I_2$, the main term does not appear from the integral by the assumption $0 \leq \alpha <1/2$ and the sum over $k$
is estimated as $ O(T^{3/4} \log^2 T)$.

Now it remains to evaluate the sum over $k$ in (\ref{Th4-I1-main-sum}).  Let
$$
S= \sum_{(\frac{T}{2\pi})^{3/2-\alpha} \leq k \leq (\frac{T}{\pi})^{3/2-\alpha}} \frac{d_3(k)}{k^{1/2}}
k^{\frac{1}{3-2\alpha}}e^{-(\frac32-\alpha)ik^{\frac{1}{3/2-\alpha}}}.
$$
By partial summation we may have
\begin{equation} \label{S-hyouka}
S \ll T^{\frac{\alpha}{2}-\frac14} \max_{(\frac{T}{2\pi})^{3/2-\alpha} \leq T' \leq (\frac{T}{\pi})^{3/2-\alpha}}
\left| \sum_{(\frac{T}{2\pi})^{3/2-\alpha} \leq k \leq T'} d_3(k)e^{-(3/2-\alpha)ik^{\frac{1}{3/2-\alpha}}}\right|.
\end{equation}
Considering the definition of $d_3(k)$, it is reduced to the estimate of the sum of the form
$$
S_1:=\sum_{T_1 \leq k_1k_2k_3 \leq 2T_1} e^{2\pi ic(k_1k_2k_3)^{\delta}},
$$
where $\delta=\frac{1}{3/2-\alpha}$, $c$ is a real constant and
$(\frac{T}{2\pi})^{3/2-\alpha} \leq T_1 \leq \frac12(\frac{T}{\pi})^{3/2-\alpha}$.
Since $\delta \neq 0, 1 $ we can apply Lemma 3 (the theorem of Robert and Sargos).
Divide the interval $[T_1,2T_1]$ into $O(\log^3 T)$ subintervals of the form $[H, 2H]\times [N,2N]\times [M,2M]$.
By symmetry of $k_j$, we can assume that $M$ is the largest, hence $M \gg T_1^{1/3}$.
Now applying Lemma 3 to the sum $S_1$ by taking $X=(HNM)^{\delta} \asymp T_1^{\delta}$, we find that
\begin{equation} \label{S1-hyouka}
S_1 \ll T_1^{1+\varepsilon}(T_1^{(\delta-\frac43)/4}+T_1^{-1/6}+T_1^{-\delta})  \ll T_1^{2/3+\delta/4+\varepsilon}.
\end{equation}
Here the last inequality follows from the assumption $0 \leq \alpha <1/2.$
By (\ref{S-hyouka}), (\ref{S1-hyouka}) and $T_1 \asymp T^{3/2-\alpha}, \ \delta=\frac{1}{3/2-\alpha}$ we find that
$$
S \ll T^{1-\frac{\alpha}{6}+\varepsilon}.
$$
This proves the assertion in the case $0 \leq \alpha <1/2$.

In the case $-1/2< \alpha \leq 0$, we take $x=\frac12(\frac{t}{2\pi})^{3/2}$ and $y=2(\frac{t}{2\pi})^{3/2}$.
Then the main term arises from the integral corresponding $I_2$ and the assertion is proved similarly.
We omit the details in this case.

\begin{flushleft}
Xiaodong Cao \\
Department of Mathematics and Physics, \\
Beijing Institute of Petro-Chemical Technology,\\
Beijing, 102617, P. R. China \\
e-mail: caoxiaodong@bipt.edu.cn

\medskip

Yoshio Tanigawa\\
Graduate School of Mathematics,\\
Nagoya University, \\
Nagoya, 464-8602, Japan\\
e-mail: tanigawa@math.nagoya-u.ac.jp

\medskip

Wenguang Zhai\\
Department of Mathematics, \\
China University of Mining and Technology, \\
Beijing 100083, P. R. China\\
e-mail: zhaiwg@hotmail.com
\end{flushleft}


\begin{thebibliography}{99}

\bibitem{Chandra} K. Chandrasekharan, Arithmetical Functions, Springer-Verlag, New York, Berlin, Heielberg, 1970.

\bibitem{I0} A. Ivi\'c,  The Theory of the Riemann Zeta-Function, Wiley \& Sons, New York 1985 (2nd ed. Dover, Mineora, 2003)


\bibitem{I2} A. Ivi\'c,  On the integral of Hardy's function, Arch. Math. {\bf 83} (2004), 41--47.

\bibitem{I3} A. Ivi\'c,  On the divisor function and the Rimann zeta-faunction in short intervals, Ramanujan J. {\bf 19}
(2009), 207--224.

\bibitem{I4} A. Ivi\'c,  The Theory of Hardy's $Z$-Function, Cambridge University Press, Cambridge, 2013.

\bibitem{I5} A. Ivi\'c,  On a cubic moment of Hardy's function with a shift, 
Exploring the Riemann Zeta Function, eds. H.L. Montgomery, A. Nikeghbali,
 M.T. Rassias, Springer Verlag, Berlin etc., 2017, pp. 99-112.

\bibitem{IZ} A. Ivi\'c and  W. Zhai, On certain integrals involving the Dirichlet divisor problem,
arXiv:1711.09589.

\bibitem{J1} M. Jutila,  Atkinson's formula for Hardy's function, J. Number Theory {\bf 129} (2009), 2853--2878.

\bibitem{J2} M. Jutila, An asymptotic formula for the primitive of Hardy's function, Ark. Math. {\bf 49} (2011), 97--107.

\bibitem{KV} A.A. Karatsuba and S.M. Vononin,  The Riemann Zeta-Function, Walter de Gruiter, Berlin, New York, 1992.

\bibitem{K} M.A. Korolev,  On the integral of Hardy's function, Izv. Math. {\bf 72} (2008), 429--478.

\bibitem{RS} O. Robert and P. Sargos, Three-dimensional exponential sums with monomials, J. reine angew. Math. {\bf 591}
(2006), 1--20.

\end{thebibliography}
\end{document}